\documentclass{article}

\input{xy}
\xyoption{all}

\usepackage{
hyperref
,srcltx
,verbatim
,varioref
,enumerate
,amsfonts
,amsthm
,amsmath
,MnSymbol
}
\usepackage[latin1]{inputenc}

		\newcommand{\category}[1]{\mathbf{#1}}

\newcommand{\Proj}{\category{Proj}} 
\newcommand{\Free}{\category{Free}} 
	
    \newcommand{\Sets}{\category{Sets}}

\newcommand{\Mod}{\category {Mod}} 

	    \newtheoremstyle{Normal}
  {}
  {}
  {}
  {}
  {\bfseries}
  {.}
  { }
  {}

    \theoremstyle{Normal} 

    \newtheorem{Defi}{Definition}[section]

    \newtheorem{Conj}[Defi]{Conjecture}
    \newtheorem{Bem}[Defi]{Remark}
    \newtheorem{Bsp}[Defi]{Example}
     
    \newtheorem{Axio}[Defi]{Axiom}
    \newtheorem{Assu}[Defi]{Assumption}
    \newtheorem{Ques}[Defi]{Question}
    \newtheorem{Comp}[Defi]{Complements}

		\theoremstyle{remark}
	
    \theoremstyle{plain} 

    \newtheorem{Satz}[Defi]{Proposition}
    \newtheorem{DefiTheo}[Defi]{Definition and Theorem}
    \newtheorem{Theo}[Defi]{Theorem}
    \newtheorem{Folg}[Defi]{Corollary}
    \newtheorem{Lemm}[Defi]{Lemma}
    \newtheorem{DefiLemm}[Defi]{Definition and Lemma}

\newcommand{\refit}[1]{(\ref{item_#1})}

\newcommand{\conj}{\begin{Conj}} 			\newcommand{\xconj}{\end{Conj}}											    
\newcommand{\ques}{\begin{Ques}} 			\newcommand{\xques}{\end{Ques}}											    
\newcommand{\axio}{\begin{Axio}} 			\newcommand{\xaxio}{\end{Axio}}											    
\newcommand{\assu}{\begin{Assu}} 			\newcommand{\xassu}{\end{Assu}}											    \newcommand{\refas}[1]{Assumption \ref{assu_#1}}
\newcommand{\bem}{\begin{Bem}} 			\newcommand{\xbem}{\end{Bem}}											    
\newcommand{\rema}{\begin{Bem}} 			\newcommand{\xrema}{\end{Bem}}											    \newcommand{\refre}[1]{Remark \ref{rema_#1}}
\newcommand{\defi}{\begin{Defi}} 			\newcommand{\xdefi}{\end{Defi}}										\newcommand{\refde}[1]{Definition \ref{defi_#1}}
\newcommand{\defitheo}{\begin{DefiTheo}} \newcommand{\xdefitheo}{\end{DefiTheo}} 
\newcommand{\defilemm}{\begin{DefiLemm}} \newcommand{\xdefilemm}{\end{DefiLemm}} 

\newcommand{\lemm}{\begin{Lemm}}			\newcommand{\xlemm}{\end{Lemm}}											\newcommand{\refle}[1]{Lemma \ref{lemm_#1}}
\newcommand{\comp}{\begin{Comp}}			\newcommand{\xcomp}{\end{Comp}}											

\newcommand{\satz}{\begin{Satz}}			\newcommand{\xsatz}{\end{Satz}}										

\newcommand{\prop}{\begin{Satz}}			\newcommand{\xprop}{\end{Satz}}										\newcommand{\refpr}[1]{Proposition \ref{prop_#1}}
\newcommand{\theo}{\begin{Theo}}			\newcommand{\xtheo}{\end{Theo}}											\newcommand{\refth}[1]{Theorem \ref{theo_#1}}
\newcommand{\bsp}{\begin{Bsp}}				\newcommand{\xbsp}{\end{Bsp}}												
\newcommand{\exam}{\begin{Bsp}}				\newcommand{\xexam}{\end{Bsp}}												\newcommand{\refex}[1]{Example \ref{exam_#1}}
\newcommand{\folg}{\begin{Folg}}				\newcommand{\xfolg}{\end{Folg}}
				 			
\newcommand{\cor}{\begin{Folg}}				\newcommand{\xcor}{\end{Folg}}
\newcommand{\mycomment}{}

\newcommand{\eqnarra}{\begin{eqnarray}}				\newcommand{\xeqnarra}{\end{eqnarray}}
\newcommand{\eqnarr}{\begin{eqnarray*}}				\newcommand{\xeqnarr}{\end{eqnarray*}}

\newcommand{\eqn}{\begin{equation}} 		\newcommand{\xeqn}{\end{equation}}
\newcommand{\refeq}[1]{(\ref{eqn_#1})}

\newcommand{\mylabel}[1]{\label{#1}}

\newcommand{\status}[1]{}
\newcommand{\explain}[1]{\stackrel{\text{\refeq{#1}}}{=}}

\newcommand{\OF}{{\mathcal{O}_F}}

\renewcommand{\log}{\mathrm{log}\,}

\newcommand{\lr}{{\longrightarrow}}
\renewcommand{\r}{\rightarrow}

\renewcommand{\t}{{\otimes}}

\newcommand{\Q}{\mathbb{Q}}
\newcommand{\Qp}{\mathbb{Q}_p}

\newcommand{\CC}{\mathbb{C}}
\newcommand{\RR}{\mathbb{R}}
\newcommand{\Z}{\mathbb{Z}}
\newcommand{\Zh}{\widehat{\mathbb{Z}}}

\newcommand{\Zp}{\mathbb{Z}_p}
\newcommand{\Zinf}{\mathbb{Z}_\infty}
\newcommand{\Zinff}{\mathbb{Z}_{(\infty)}}
\newcommand{\Fp}{\mathbb{F}_p}
\newcommand{\Finf}{\mathbb{F}_\infty}
\newcommand{\Fone}{\mathbb{F}_1}
\newcommand{\Fzero}{\mathbb{F}_0}
\newcommand{\pp}{\mathfrak{p}}

\renewcommand{\H}{\mathrm{H}}

\newcommand{\x}{{\times}}

\newcommand{\Beweis}{{\normalfont} \textbf{Proof}}

\newcommand{\lcs}{\textit{loc.~cit.}\ }

\newcommand{\Spec}{\mathrm{Spec} \; \! }

\newcommand{\SpecZ}{{\Spec} {\Z}}

\newcommand{\SpecZh}{{\Spec} {\Zh}}

\newcommand{\rk}{\operatorname{rk}} 

\newcommand{\weq}{\stackrel{\sim}\r}
\newcommand{\cofi}{\rightarrowtail}
\newcommand{\fib}{\twoheadrightarrow}
\newcommand{\fibl}{\twoheadleftarrow}
\newcommand{\sgn}{\operatorname{sgn}} 

\newcommand{\pist}[1]{\pi^\mathrm{s}_{#1}}
\newcommand{\Hom}{\mathrm{Hom}}
\newcommand{\End}{\mathrm{End}}

\newcommand{\id}{\mathrm{id}}

\newcommand{\im}{\operatorname{im}} 
\newcommand{\ab}{\mathrm{ab}} 

\newcommand{\Aut}{\mathrm{Aut}}

\newcommand{\coker}{\operatorname{coker}}

\newcommand{\Id}{\mathrm{Id}}

\newcommand{\loc}[2]{[}

\newcommand{\pr}{\begin{proof}[\Beweis: ]}
\newcommand{\pf}{\pr}
\newcommand{\xpf}{\end{proof}}

\theoremstyle{definition}

\theoremstyle{definition}

\theoremstyle{definition}

\theoremstyle{definition}

\theoremstyle{definition}

\theoremstyle{definition}

\theoremstyle{definition}

\theoremstyle{definition}

\theoremstyle{definition}

\theoremstyle{definition}

\theoremstyle{definition}

\theoremstyle{definition}

\theoremstyle{definition}

\theoremstyle{definition}

\newcommand{\cO}{\mathcal{O}}

\newcommand{\cC}{\mathcal{C}}

\newcommand{\GL}{\mathrm{GL}}

\renewcommand{\mycomment}{\tiny}

\usepackage{titlesec}
\usepackage[]{titletoc}
 
\titlecontents{subsubsubsection}[9em]{}{\contentslabel{3.9em}}%
{\hspace*{-1.2em}}{\titlerule*[0.675pc]{.}\contentspage}
 
\makeatletter
\newcounter{subsubsubsection}[subsubsection]
\renewcommand{\thesubsubsubsection}{\thesubsubsection.\@arabic\c@subsubsubsection}
 
\titleclass{\subsubsubsection}{straight}[\subsubsection]
\titleformat{\subsubsubsection}{\bf}{\thetitle}{1em}{}[]						
\titlespacing{\subsubsubsection}{0pt}{3.25ex plus 1ex minus 0.2ex}{1.5ex plus 0.2ex} 
 
\makeatother

\begin{document}

\author{Jakob Scholbach\footnote{Universit{\"a}t M{\"u}nster, Mathematisches Institut, Einsteinstr. 62, D-48149 M{\"u}nster, Germany.}
}

\title{Algebraic $K$-theory of the infinite place}


\maketitle


\begin{abstract}
We show that the algebraic $K$-theory of generalized archimedean valuation rings occurring in Durov's compactification of the spectrum of a number ring is given by stable homotopy groups of certain classifying spaces. We also show that the ``residue field at infinity'' is badly behaved from a $K$-theoretic point of view. 
\end{abstract}

\section{Introduction}
In number theory, it is a universal principle that the spectrum of $\Z$ should be completed with an infinite prime. This is corroborated, for example, by Ostrowski's theorem, the product formula
$$\prod_{p \leq \infty} |x|_p = 1, \ \ x \in \Q^\x,$$
the Hasse principle, Artin-Verdier duality, and functional equations of $L$-functions.

This ``compactification'' $\SpecZh := \SpecZ \cup \{ \infty \}$ was just a philosophical device until recently: Durov has proposed a rigorous framework which allows for a discussion of, say, $\Zinff$, the local ring of $\SpecZh$ at $p=\infty$ \cite{Durov}. The purpose of this work is to study the $K$-theory of the so-called generalized rings intervening at the infinite place. 

Algebraic $K$-theory is a well-established, if difficult, invariant of arithmetical schemes. For example, the pole orders of the Dedekind $\zeta$-function $\zeta_F(s)$ of a number field $F$ are expressible by the ranks of the $K$-theory groups of $\OF$, the ring of integers. By definition, $K$-theory only depends on the category of projective modules over a ring. Therefore, this interacts nicely with Durov's theory of \emph{generalized rings} which describes (actually: defines) such a ring $R$ by defining its free modules. For example, the free $\Zinff$-module of rank $n$ is defined as the $n$-dimensional octahedron, i.e.,
$$\Zinff(n) := \{ (x_1, \dots, x_n) \in \Q^n, \sum_i |x_i| \leq 1 \}.$$
The abstract theory of such modules is a priori more complicated than in the classical case since $\Zinff$-modules fail to build an abelian category. Nonetheless, using Waldhausen's $S_\bullet$-construction it is possible to study the \emph{algebraic $K$-theory} of $\Zinff$ and similar rings occuring for other number fields (\refth{Waldhausen}, \refde{K}). 

\hspace{0cm}\\
\noindent \textbf {\refth{Khigher}.}
\textit{ 
The $K$-groups of $\Zinff$ are given by
$$K_i(\Zinff) = \pi_{i}^\mathrm{s} (B \mu_2 \sqcup \{*\}, *) = 
\begin{cases}
\Z & i = 0 \ \text{(Durov \cite[10.4.19]{Durov})} \\
\Z / 2 \oplus \mu_2 & i=1 \\
\mathrm{a\ finite\ group} & i > 1.
                                             \end{cases}$$
The $\Z/2$-part in $K_1$ stems from the first stable homotopy group $\pi_1^s$, while $\mu_2 = \{\pm 1 \}$ arises as the subgroup of $\Zinff$ of elements of norm $1$, i.e., the subgroup of (multiplicative) units of $\Zinff$. The finite $K$-group for $i>1$ is the abutment of an Atiyah-Hirzebruch spectral sequence.
}
\hspace{0cm}\\

This theorem is proven for more general generalized valuation rings including $\OF_{(\sigma)}$, the ring corresponding to an infinite place $\sigma$ of a number field $F$. In this case the group $\mu_2$ above is replaced by the group $\{x \in F, |\sigma(x)| = 1\}$. The basic point is this: the only admissible monomorphisms (i.e., the ones occurring in the $S_\bullet$-construction of $K$-theory) 
$$\Zinff(1) = [-1, 1] \cap \Q \r \Zinff(2)$$
are given by mapping the interval to one of the two diagonals of the lozenge. Thereby, the Waldhausen category structure on free $\Zinff$-modules turns out to be equivalent to the one of finitely generated pointed $\{\pm1\}$-sets, whose $K$-theory is well-known. In the course of the proof we also show that other plausible definitions, such as the $S^{-1} S$-construction, the $Q$-construction, and the $+$-construction yield the same $K$-groups. 

We finish this note by pointing out two $K$-theoretic differences of the infinite place: we show that $K_0(\Finf) = 0$ (\refpr{K0Finf}), as opposed to $K_0(\Fp) = \Z$. Also, the completions at infinity are not well-behaved from a $K$-theoretic viewpoint. These remarks raise the question whether the ``local'' ring $\Zinff$ should be considered regular or, more precisely, whether  
$$K_0(\Zinff) \r K'_0(\Zinff) := \Z[\text{finitely presented }\Zinff-\Mod] / \text{short exact sequences} \ \ $$
is an isomorphism. Unlike in the classical case, there does not seem to be an easy resolution argument in the context of Waldhausen categories. Another natural question is whether there is a Mayer-Vietoris sequence of the form
$$K_i(\Zh) \r K_i(\Z) \oplus K_i(\Zinff) \r K_i(\Q) \r K_{i-1}(\Zh),$$
where $\Zh$ is a generalized scheme obtained by glueing $\SpecZ$ and $\Spec \Zinff$ along $\Spec \Q$. 
The usual proof of this sequence proceeds by the localization sequence, which is not available in our context.

Throughout the paper, we use the following \emph{notation}\mylabel{notation}: $F$ is a number field with ring of integers $\OF$. Finite primes of $\OF$ are denoted by $\pp$. We write $\Sigma_F$ for the set of real and pairs of complex embeddings of $F$. The letter $\sigma$ usually denotes an element of $\Sigma_F$. It is referred to as an infinite prime of $\OF$.

\emph{Acknowledgement:} 
I would like to thank Fabian Hebestreit for a few helpful discussions. I also thank the referee for suggesting a number of improvements.

\section{Generalized rings}
In a few brushstrokes, we recall the definition of generalized rings and their modules and some basic properties. Everything in this section is due to Durov. All references in brackets refer to \cite{Durov}, where a much more detailed discussion is found. 

A monad in the category of sets is a functor $R : \Sets \r \Sets$ together with natural transformations $\mu: R \circ R \r R$ and $\epsilon : \Id \r R$ required to satisfy an associativity and unitality axiom akin to the case of monoids. We will write $R(n) := R(\{1, \dots, n\})$. An $R$-module is a set $X$ together with a morphism of monads $R \r \End(X)$, where the endomorphism monad $\End(X)$ satisfies $\End(X)(n) = \Hom_{\Sets}(X^n, X)$. In other words, $X$ is endowed with an action 
$$R(n) \x X^n \r X$$
satisfying the usual associativity conditions. Thus, $R(n)$ can be thought of as the $n$-ary operations (acting on any $R$-module).

\defi (Durov [5.1.6])
A \emph{generalized ring} is a monad $R$ in the category of sets satisfying two additional properties:
\begin{itemize}
\item $R$ is \emph{algebraic}, i.e., it commutes with filtered colimits. Since every set is the filtered colimit of its finite subsets, this implies that $R$ is determined by $R(n)$ for $n \geq 0$ [4.1.3].
\item $R$ is \emph{commutative}, i.e., for any $t \in R(n)$, $t' \in R(n')$, any $R$-module $X$ (it suffices to take $X = R(n \x n')$) and $A \in X^{n\x n'}$, we have 
$$t(t'(A)) = t'(t(A)),$$ 
where on the left hand side $t'(A) \in X^n$ is obtained by letting act $t'$ on all rows of $A$ and similarly (with columns) on the right hand side. 
\end{itemize}
\xdefi

For a unital associative ring $R$ (in the sense of usual abstract algebra), let 
$$R(S) := \oplus_{s \in S} R$$
be the free $R$-module of rank $\sharp S$, where $S$ is any set. The addition and multiplication on $R$ turn this into an (algebraic) monad which is commutative iff $R = R(1)$ is [3.4.8]. Indeed, the required map 
\eqn \mylabel{eqn_R1}
R(1) \x R(1) \r R(1)
\xeqn
is just the multiplication in $R$, while the addition is reformulated as 
$$R(2) \x (R(1) \x R(1)) \r R(1), ((x_1, x_2), (y_1, y_2)) \mapsto \sum x_i y_i.$$
Note that \refeq{R1} is required to exist for any monad, so multiplication is in a sense more fundamental than addition, which requires the particular element $(1, 1) \in R(2)$ [3.4.9]. 

Reinterpreting a ring as a monad in this way defines a functor from commutative rings to generalized rings, which is easily seen to be fully faithful: given two classical rings $R$, $R'$, and a map of monads, i.e., a collection of maps $R(n) = R^n \r R'(n) = R'^n$, one checks that the maps for $n \geq 2$ are determined by $R \r R'$. In the same vein, $R$-modules in the classical sense are equivalent to $R$-modules (in the generalized sense). Henceforth, we will therefore not distinguish between classical commutative rings and their associated generalized rings.  

The initial generalized ring is the monad $\Fzero: \Sets \r \Sets$, $M \mapsto M$. Its modules are just the same as sets. The monad $\Sets \ni M \mapsto M \sqcup \{ * \}$ is denoted $\Fone$. Neither of these two generalized rings is induced by a classical ring. See \refde{cO} for our main example of a non-classical ring.    

Given a morphism $\phi: R \r S$ of generalized rings, the forgetful functor $\Mod(S) \r \Mod(R)$ between the module categories has a left adjoint $\phi^*: \Mod(R) \r \Mod(S)$ called \emph{base change}. We also denote it by $-\t_R S$. Being a left adjoint, this functor preserves colimits [4.6.19].\mylabel{colimits} For example, for a generalized ring $R$, the unique map $\Fzero \r R$ of generalized rings induces an adjunction
$$\Sets = \Mod(\Fzero) \leftrightarrows \Mod(R): \text{forget}$$
Its left adjoint is explicitly given by $X \mapsto R(X)$, the so-called \emph{free $R$-module} on some set $X$. That is, 
$$\Hom_{\Mod(R)}(R(X), M) = \Hom_{\Sets}(X, M),$$
as in the classical case.  

Coequalizers and arbitrary coproducts exist in $\Mod(R)$, for any generalized ring $R$ [4.6.17]. Therefore, arbitrary colimits exist. Base change functors $\phi^*$ commute with coequalizers. 
Moreover, arbitrary limits exist in $\Mod(R)$, and commute with the forgetful functor $\Mod(R) \r \Sets$ [4.6.1]. 

An $R$-module $M$ is called \emph{finitely generated} if there is a surjection $R(n) \fib M$ for some $0 \leq n < \infty$  [4.6.9]. 
Unless the contrary is explicitly mentioned, all our modules are supposed to be finitely generated over the ground generalized ring in question. 
An $R$-module $M$ is \emph{projective} iff it is a retract of a free module, i.e., if there are maps $M \stackrel i \r R(n) \stackrel p \r M$ 
with $p i = \id_M$. As in the classical case this is equivalent to the property that for any surjection of $R$-modules $N \fib N'$, $\Hom_{\Mod(R)}(M, N)$ maps onto $\Hom_{\Mod(R)}(M, N')$  [4.6.23]. The categories of (finitely generated) free and projective $R$-modules are denoted $\Free(R)$ and $\Proj(R)$, respectively. 

As usual, an \emph{ideal} $I$ of $R$ is a submodule of $R(1)$. A proper ideal $I \subsetneq R(1)$ is called \emph{prime} if $R(1) \backslash I$ is multiplicatively closed [6.2.2]. 

\section{Archimedean valuation rings}
\subsection{Definitions}
Let $K$ be an integral domain equipped with a norm $|-|: K \r \RR^{\geq 0}$. We will write $Q$ for the quotient field of $K$. We put $E := \{ x \in K, |x| = 1 \}$. We also write $|x|$ for the $L^1$-norm on $K^n$, i.e., $|x| = \sum_i |x_i|$. Throughout, we assume:
\assu \mylabel{assu_basic}
\begin{enumerate}[(A)]
\item \mylabel{item_assA} $|K^\x| = \{|k|, k \in K^\x\} \subset \RR^{\geq 0}$ is dense. 
\item \mylabel{item_assB} $E \subset K^\x$.
\end{enumerate}
\xassu

\defi \mylabel{defi_cO} 
The \emph{(generalized) valuation ring} associated to $(K, |-|)$ is the submonad $\cO$ of $K$ given by
$$\cO(S) := \left \{ x = (x_s) \in \bigoplus_{s \in S} K, |x| := \sum_{s \in S} |x_s| \leq 1 \right \}.$$
\xdefi
This is clearly algebraic. Moreover, the multiplication of the monad, i.e., $\cO \circ \cO \r \cO$ is well-defined by restricting the one of $K$ (and is therefore commutative):
$$\cO(\cO(n)) = \left \{ (y_x) \in \bigoplus_{x \in \cO(n)} K, \sum_x |y_x| \leq 1 \right \} \r \cO(n)$$
sends $(y_x)$ to (the finite sum) $\sum_x y_x \cdot x$. A priori, this expression is an element of $K^n$, only, but is actually contained in $\cO(n)$ since
$$\left |\sum_x y_x \cdot x \right | \leq \left (\sum_x |y_x| \right ) \cdot \sup |x| \leq 1.$$
In the case of an archimedean valuation, this definition of $\cO$ is the one of Durov \cite[5.7.13]{Durov}. For non-archimedean valuations, Durov's original definition gives back the (generalized ring corresponding to the) ordinary ring $\{x \in K, |x| \leq 1\}$ which is different from \refde{cO} (see \refex{1}).   

By definition, an $\cO$-module $M$ is therefore a set such that an expression $\sum_{i=1}^n \lambda_i m_i$ is defined for $n \geq 0$, $m_i \in M$, $\lambda_i \in K$ such that $\sum |\lambda_i| \leq 1$, obeying the usual laws of commutativity, associativity and distributivity. Maps $f: M \r N$ of $\cO$-modules are described similarly: they satisfy $f(\sum_i \lambda_i m_i) = \sum_i \lambda_i f(m_i)$. The set $\{0\}$, with its obvious $\cO$-module structure is both an initial and terminal $\cO$-module. Given a map $f: M' \r M$ of $\cO$-modules, the (co)kernel is defined to be the (co)equalizer of the two morphisms $f$ and $M' \r 0 \r M$. As was noted above, the forgetful functor $\cO-\Mod \r \Sets$ preserves limits, so the kernel $\ker f$ is just $f^{-1}(0)$. The cokernel is described by the following proposition. Also see \refre{split} for an explicit example of a cokernel computation.

\prop \mylabel{prop_cokernel}
Given a map $f: M' \r M$ of $\cO$-modules, the cokernel is given by    
\eqn
\mylabel{eqn_cokernel}
\coker (f) = M \left / \sim \right .,
\xeqn
where $\sim$ is the equivalence relation generated by $\sum_{i \in I} \lambda_i m_i \sim \sum_{i \in I} \lambda_i \tilde m_i$, where $I$ is any finite set, $\lambda = (\lambda_i) \in \cO(\sharp I)$ and $m_i, \tilde m_i \in M$ are such that either $m_i = \tilde m_i$ or both $m_i$, $\tilde m_i \in f(M') \subset M$. This set is endowed with the $\cO$-action via the natural projection $\pi: M \r \coker (f)$. 
\xprop
\pf
This follows from the description of cokernels given in \cite[4.6.13]{Durov}. It is also easy to check the universal property directly: 
we clearly have $\pi \circ f = 0$. Given a map $t: M \r T$ of $\cO$-modules such that $t f = 0$, we need to see that $t$ factors uniquely through $\coker f$. The unicity of the factorization is clear since $M \r \coker f$ is onto. The existence is equivalent to $t(m_1) = t(m_2)$ whenever $\pi (m_1) = \pi (m_2)$. This is obvious from the definition of the equivalence relation $\sim$ above. 
\xpf

The base change functor resulting from the monomorphism $\cO \subset K$ of generalized rings is denoted
$$(-)_K : \Mod(\cO) \r \Mod(K).$$
Actually, using \refas{basic}, we may pick $t \in K^\x$ such that $|t|<1$. Then, $K$ is the unary localization $K = \cO[1/t]$. This is shown in \cite[6.1.23]{Durov} for $K = \RR$. The proof for a general domain is the same. Therefore $K$ is flat over $\cO$, so $(-)_K$ preserves finite limits, in particular kernels \mylabel{localization}\cite[6.1.2, 6.1.8]{Durov}. Recall from p.\ \pageref{colimits} that $(-)_K$ also preserves colimits, such as cokernels. 

Let $E(n) := \{ x \in K(n) = K^n, |x| = 1 \}$ be the ``boundary'' of $\cO(n)$. (This is merely a collection of sets, not a monad.) We write $\cO$ for $\cO(1)$ and $E$ for $E(1)$, if no confusion arises. In particular, $x \in \cO$ means $x \in \cO(1)$. The $i$-th standard coordinate vector $e_i = (0, ..., 1, ..., 0)$ is called a \emph{basis vector} of $\cO(n)$ ($1 \leq i \leq n$). 

\exam \mylabel{exam_1}
Let $F$ be a number field with ring of integers $\OF$. We fix a complex embedding $\sigma: F \r \CC$ and take the norm $|-|$ induced by $\sigma$. Let $K$ be either $\cO_F[1/N]$ where $N \in \Z$ has at least two distinct prime divisors, or $F$, or $\widehat F^{\sigma}$, the completion of $F$ with respect to $\sigma$. The respective generalized valuation rings will be denoted $\cO_{F, 1/N, (\sigma)}$, $\cO_{F, (\sigma)}$, and $\cO_{F,\sigma}$, respectively. For example, $\cO_{F, (\sigma)} = \cO_{F, (\overline \sigma)}$. \refas{basic}\refit{assA} is satisfied: for $\cO_F[1/N]$, pick two distinct prime divisors $p_1 \neq p_2$ of $N$. The elements $p_1^{n_1} p_2^{n_2} \in K$ are invertible for any $n_1$, $n_2 \in \Z$. The subgroup $\{ \log (| p_1^{n_1} p_2^{n_2} |), n_i \in \Z \} \subset \RR$ is dense: otherwise it was cyclic, in contradiction to the $\Q$-linear independence of $\log p_1$ and $\log p_2$ (Gelfand's theorem).

As for \refas{basic}\refit{assB}, let $x \in K$ with $|x|=1$. If $\sigma$ is a real embedding, $x = \pm |x| = \pm 1$. If $\sigma$ is a complex embedding, let $\overline \sigma$ be its complex conjugate and $\overline x \in K$ be such that $\sigma (\overline x)=\overline \sigma(x)$. Then $\sigma(x) \sigma(\overline x)=\sigma(x) \overline{\sigma(x)} = |\sigma(x)|^2 = 1$ implies $x \in K^\x$.  

According to Durov, $\cO_{F, (\sigma)}$ is the replacement for infinite places of the local rings $\OF_{(\pp)}$ at finite places. 
However, the analogy is relatively loose, as is shown by the following two remarks: first, for $p < \infty$, let $|x|_p := p^{-v_p(x)}$ for $x \in \Q^\x$. Then the generalized ring $\Z_{|-|_p}$ (in the sense of \refde{cO}) maps injectively to the localization $\Z_{(p)}$ of $\Z$ at the prime ideal $p$, but the map is a bijection only in degrees $\leq p$. (Less importantly, \refas{basic}\refit{assA} is not satisfied for $\Z_{|-|_p}$.) 

Secondly, recall that the semilocalization $\OF_{(\pp_1, \pp_2)} = \OF_{(\pp_1)} \cap \OF_{(\pp_2)}$ at two finite primes is one-dimensional. In analogy, pick two $\sigma_1, \sigma_2 \in \Sigma_F$ and consider $\cO := \cO_{(\sigma_1)} \cap \cO_{(\sigma_2)} \subset F$, i.e., 
$$\cO(n) := \{(x_1, \dots, x_n) \in F^n, \sum_k |\sigma_i(x_k)| \leq 1 \ \ \text{ for }i=1, 2\}.$$
Let $\mathfrak p_i = \{ x \in \cO, |\sigma_i(x)| < 1 \}$ and $\mathfrak p := \{ x \in \cO, |\sigma_1(x) \sigma_2(x)| < 1 \}$. These are ideals: for example, for $x = (x_j) \in \cO(n)$, $s_1, \dots, s_n \in \mathfrak p$, we need to check $\sum s_j x_j \in \mathfrak p$: if, say, $|\sigma_1(s_1)| < 1$ then 
$$\left |\sigma_1 (\sum_j s_j x_j) \right| \leq \sum |\sigma_1 (s_j)| |\sigma_1(x_j)| < \sum |\sigma_1(x_j)| \leq 1.$$
The complement $\cO \backslash \mathfrak p = \{ x, |\sigma_1(x)| = |\sigma_2(x)| = 1 \}$ is multiplicatively closed (and contains $1$). We get a chain of prime ideals
$$0 \subsetneq \mathfrak p_1 \subset \mathfrak p \subsetneq \cO.$$
The middle inclusion is, in general, strict, namely when $F = \Q[t] / p(t)$ with some irreducible polynomial $p(t)$ having zeros $a_1, a_2 \in \CC$ with $|a_1| = 1$, $|a_2| < 1$. That is, $\Spec \cO$ is not one-dimensional. 
\xexam

\subsection{Projective and free $\cO$-modules}
In this section we gather a few facts about projective and free $\cO$-modules. We begin with a handy criterion for monomorphisms of certain $\cO$-modules (\refle{morphismscO}). \refle{particular} concerns a particular unicity property of the basis vectors $e_i = (0, \dots, 0, 1, 0, \dots, 0) \in \cO(n)$. This is used to prove \refth{projfree}: every projective $\cO$-module is free, provided that the norm is archimedean. This improves a result of Durov which treats only the cases where $\cO$ is either the ``unclompeted local ring'' of a number ring at an infinite place $\sigma$, $\cO_{F, (\sigma)}$, in the case where $\sigma$ is a real embedding or the ``completed local ring'' $\cO_{F, \sigma}$ for both real and complex places. Therefore, we only study the $K$-theory of free $\cO$-modules in this paper (but see \refre{cofinality}). We also use \refle{particular} to establish a highly combinatorial flavor of automorphisms of free $\cO$-modules (\refpr{morphismscO2}), which will later give rise to the 
computation of higher $K$-theory of $\cO$.     



\lemm \mylabel{lemm_morphismscO} (compare \cite[2.8.3.]{Durov}) 
Let $f : M' \r M$ be a map of $\cO$-modules. We suppose both $M'$ and $M$ are submodules of free $\cO$-modules. (For example, they might be projective.) Then the following are equivalent:
\begin{enumerate}[a)]
\item \mylabel{item_l2d} $f_Q: M'_Q \r M_Q$ is injective, where $Q$ is the quotient field of $K$,
\item \mylabel{item_l2c} $f_K: M'_K \r M_K$ is injective,
\item \mylabel{item_l2b} $f$ is injective (as a map of sets),
\item \mylabel{item_l2a} $f$ is a monomorphism of $\cO$-modules,
\end{enumerate}
\xlemm
\pf
Consider the diagram
$$\xymatrix{
M' \ar[d]^f \ar@{^{(}->}[r] & M'_K \ar[d]^{f_K} \ar@{^{(}->}[r] & M'_Q \ar[d]^{f_Q} \\
M  \ar@{^{(}->}[r] & M_K \ar@{^{(}->}[r] & M'_Q.
}$$
Its horizontal maps are injective since both modules are submodules of free modules and, for these, $\cO(n) \subset K(n) = K^n \subset Q(n) = Q^n$. This shows \refit{l2d} $\Rightarrow$ \refit{l2c} $\Rightarrow$ \refit{l2b}. \refit{l2b} implies \refit{l2a} since the forgetful functor $\Mod(\cO) \r \Sets$ is faithful. 
\refit{l2a} $\Rightarrow$ \refit{l2c}: by \refas{basic}, we may pick $t \in K^\x$ with $|t|<1$. Any two element of $M'_K$ are of the form $m'_1 / t^n$, $m'_2 / t^n$, where $m'_1, m'_2 \in M'$ and $n \geq 0$. Suppose that $f_K(m'_1 / t^n) = f(m'_1) / t^n$ agrees with $f_K(m'_2/t^n)$. The multiplication with $t^{-n}$ is injective on $M'_K$, since $M'$ ($M'_K$) is a submodule of a free $\cO$- ($K$-, respectively) module. Thus $f(m'_1) = f(m'_2)$ so the assumption \refit{l2a} implies our claim. 
Finally \refit{l2c} $\Rightarrow$ \refit{l2d} follows from the flatness of $Q$ over $K$. 
\xpf

The following lemma can be paraphrased by saying that the basis vectors $e_i = (0, \dots, 1, \dots 0) \in \cO(n)$ cannot be generated as a nontrivial $\cO$-linear combination of other elements of $\cO(n)$.

\lemm \mylabel{lemm_particular}
Suppose that $K$ is a field (as opposed to a domain). Suppose further that
\eqn \mylabel{eqn_special}
e_i = \sum_{j=1}^m \lambda_j f_j
\xeqn
with $f_j \in \cO(n)$ and $(\lambda_j)_j \in \cO(m)$, $\lambda_j \neq 0$. Then for each $j$, $f_j = \mu_j \cdot e_i$ with $\mu_j \in E$. 
\xlemm
\pf
The proof proceeds by induction on $m$, the case $m=1$ being trivial.

Each $f_j$ can be written as $f_j = \sum_{l=1}^n \kappa_{jl} e_l$ with $(\kappa_{jl})_l \in \cO(n)$. We get
\eqn
\mylabel{eqn_pf1}
1 = |e_i| \explain{special} |\sum \lambda_j f_j| \leq \sum |\lambda_j| |f_j| \leq \sum |\lambda_j| \leq 1.
\xeqn
Therefore equality holds throughout. 
%
%
We have $e_i = \sum_{j,l} \lambda_j \kappa_{jl} e_l$. This $K$-linear relation between the basis vectors of $K^n$ yields $1 = \sum_j \lambda_j \kappa_{ji}$. Hence
$$1 \leq \sum_j |\lambda_j \kappa_{ji}| \leq \underbrace{(\sum |\lambda_j|)}_{\explain{pf1}1} \cdot \max_j |\kappa_{ji}|.$$
On the other hand, $|\kappa_{ji}| \leq 1$, so there is some $j_0$ such that $|\kappa_{j_0i}|=1$.  
Using $\sum_l |\kappa_{j_0l}| \leq 1$ we see $\kappa_{j_0l} = 0$ for all $l\neq i$, thus $f_{j_0} = \kappa_{j_0i} e_i.$ Put $\mu_{j_0} := \kappa_{j_0i} (\in E)$, so  
$$(1-\lambda_{j_0} \mu_{j_0}) e_i = \sum_{j \neq j_0} \lambda_{j} f_j$$
holds. If $|\lambda_{j_0} \mu_{j_0}| = 1$, we are done since all other $\lambda_j$, $j \neq j_0$ must vanish in this case. If $|\lambda_{j_0} \mu_{j_0}| < 1$, then
$$e_i = \sum_{j \neq j_0} \frac{\lambda_j}{1 - \lambda_{j_0} \mu_{j_0}}	f_j.$$
This finishes the induction step since the right hand side is actually an $\cO$-linear combination of the $f_j$, for
$$\sum_{j \neq j_0} |\lambda_{j} | \explain{pf1} 1 - |\lambda_{j_0}| =  1 - |\lambda_{j_0} \mu_{j_0}| \leq |1 - \lambda_{j_0} \mu_{j_0}|.$$
\xpf

\theo \mylabel{theo_projfree} 
Suppose that the norm $|-|$ giving rise to the generalized valuation ring $\cO$ is archimedean. Then every projective $\cO$-module $M$ is free.  
\xtheo
\pf
Let $K'$ be the completion (with respect to the norm $|-|$) of $Q$, the quotient field of $K$. By Ostrowski's theorem, we have either $K' = \RR$ or $K' = \CC$ (with their usual norms). Let us write $-' := - \t_\cO \cO'$, where $\cO' := \cO_{K'}$ is the generalized valuation ring belonging to $K'$. We consider the following maps of $\cO'$-modules, where $O_i$ are certain free $\cO$-modules that are defined in the course of the proof:
$$O_3' \r O_2' \r O_1' \stackrel{p'}\lr M' \stackrel{\phi, \cong} \lr O_0'.$$
First, $M'$ is a projective $\cO'$-module: given a projector $p : O_1 := \cO(n_1) \r \cO(n_1)$ with $M = \im p$, we get $M' = \im p'$. By the afore-mentioned result of Durov \cite[10.4.2]{Durov}, there is an isomorphism of $\cO'$-modules, $\phi: M' \stackrel \cong \r O_0' := \cO'(n_0)$. The composition $\phi \circ p'$ is surjective, so for any basis vector $e_i \in O_0'$ ($1 \leq i \leq n_0$), there is some $\cO'$-linear combination $\sum_{j \leq n_1} \lambda_{ij} e_j$ mapping to $e_i$ under $\phi p'$. Thus, $\sum_j \lambda_{ij} \phi p' (e_j) = e_i$. Therefore, by \refle{particular}, $\phi p' (e_j) \in E' \cdot e_i$ for each $j$. 
Here $E' = \{ x \in \cO', |x| = 1\}$ (which is $S^1 \subset \CC$ or $\{\pm 1\} \subset \RR$ depending on $K'$). 
We put $O_2 := \sqcup_{j_2 \in J_2} e_{j_2} \cO = \cO(J_2)$, where the coproduct runs over 
$$J_2 := \{1 \leq j_2 \leq n_1, \phi p'(e_{j_2}) \in E' e_i \ \ \text {for some }i \leq n_0\}.$$ 
The inclusion $J_2 \subset \{1, \dots, n_1\}$ induces a ($\cO$-linear!) injection $f_{21}: O_2 \r O_1$. According to the previous remark, $O'_2 \stackrel{\phi p' f_{21}'} \lr O'_1$ is surjective. Consider the map $J_2 \r \{1, \dots, n_0 \}$ which maps $j_2$ to the (unique) $i$ with $e_i \in E' \phi p'(e_{j_2})$. This map is onto. By \refas{basic}, we may pick some $J_3 \subset J_2$ on which it is a bijection. Let $f_{32}: O_3 := \sqcup_{j_3 \in J_3} e_{j_3} \cO = \cO(J_3) \r O_2 = \cO(J_2)$ be the map induced by $J_3 \subset J_2$. Set $f_{31} = f_{21} \circ f_{32}$. Then the composition $O'_3 \stackrel{f_{31}'} \cofi O'_1 \stackrel{p'}\r M' \stackrel {\phi, \cong} \lr O'_0$ is an isomorphism of $\cO'$-modules. Note that $f_{31}$ and $p$ are $\cO$-linear maps, but $\phi$ is defined over $\cO'$, only. Writing $v := p \circ f_{31}$, we must show the implication
$$v' \text{ isomorphism} \Rightarrow v \text { isomorphism}.$$ 

The elements $m_j := p(e_j) \in M$, $j \leq n_1$, generate $M$. The map $v' \t_{\cO'} K' = v_{Q} \t_{Q} K'$ is an isomorphism of $K'$-vector spaces. The inclusion of the quotient field $Q \r K'$ is fully faithful, so that $v_{Q}$ is also an isomorphism. Hence there is some $k_j = a_j / b_j \in Q \backslash \{ 0 \}$ such that $k_j m_j \in \im v$. According to \refas{basic}, we can pick some $N \in K^\x$ such that $|a_j / N|, |b_j / N| \leq 1$ for all $j$. Then $m_j a_j / N \in \im v$. Similarly, pick some $t \in \cO$ with $0 < |t| \leq \min_j |a_j / N|$. Then $t M \subset \im v$. 

To show the surjectivity of $v$, we fix $m \in M$ and pick some $o_3 \in O_3$ with $tm = v(o_3)$. Since $M \subset M'$ and $v'$ is an isomorphism, there is a unique $\tilde o'_3 \in O'_3$ with $v'(\tilde o'_3) = m$. Hence $v(o_3) = v'(o_3) = v'(t\tilde o'_3)$, so that $t \tilde o'_3 = o_3$. In other words, $o'_3 = t^{-1} o_3 \in O'_3 \cap (O_3)_K = O_3$. This shows the surjectivity of $v$. The injectivity of $v$ is clear, since $O_3 \subset O'_3$ and $v'$ is injective. Consequently, $v$ is an isomorphism.        
\xpf

\defi \mylabel{defi_Waldhausen} 
Recall that $\Free(\cO)$ is the category of (finitely generated) free $\cO$-modules. In $\Free(\cO)$ let \emph{cofibrations} ($\cofi$) be the monomorphisms whose cokernel (in the category of all $\cO$-modules) lies in $\Free(\cO)$. Morphisms which are obtained as cokernels of cofibrations are called \emph{fibrations} and denoted $\fib$. Let \emph{weak equivalences} $\weq$ be the isomorphisms. 
\xdefi

\prop \mylabel{prop_morphismscO2}
Let $f: M' \r M$ be a monomorphism of free $\cO$-modules with projective cokernel $M''$ (for example, a cofibration). Then there is a \emph{unique} isomorphism $\phi: M \cong M' \sqcup M''$ such that the following diagram is commutative
\eqn \mylabel{eqn_short.exact}
\xymatrix{
M' \ar@{>->}[r]^f \ar@{=}[d] & M \ar@{->>}[r]^\pi \ar@{.>}[d]^\phi & M'' \ar@{=}[d] \\ 
M' \ar@{>->}[r]^{\text{incl}} & M' \sqcup M'' \ar@{->>}[r]^{\text{proj}} & M''.
}.
\xeqn
\xprop

\pf 
Let $M' = \cO(n')$, $M = \cO(n)$ and let $f_i := f(e_i) \in M$, $1 \leq i \leq n'$ be the images of the basis vectors. 

We claim that 
$f$ factors through $\sqcup_{i \leq n, e_i \in f(M')} e_i \cO = \cO(\tilde n') \subset M = \cO(n)$,
where $\tilde n' := \sharp \{i \leq n, e_i \in f(M') \}$. To show this, write $f(M') \ni m' = \sum_{i \in I} \lambda_i e_i$, where all $\lambda_i \neq 0$ and the $e_i$ are the basis vectors of $M$. Put
$$m' = \underbrace{\sum_{e_i \notin f(M')} \lambda_i e_i}_{=:m'_1} + \underbrace{\sum_{e_i \in f(M')} \lambda_i e_i}_{=:m'_2}.$$
By \refas{basic}, we can pick some $t \in K^\x$ such that $|t| \leq 1/2$. Then $t m'_1 = t m' - t m'_2 \in f(M')$. Let $i$ be such that $e_i \notin f(M')$. We need to see $\lambda_i = 0$. 

We write $(-)_Q$ for the functor $- \t_\cO \cO_Q$, where $\cO_Q$ is the generalized valuation ring associated to the unique extension of the norm $|-|$ in $K$ to the quotient field $Q$ of $K$. The functor $(-)_Q$ preserves colimits, in particular $\coker (f_Q) = (\coker f)_Q$. In addition, $f_Q$ is a monomorphism by \refle{morphismscO}. The assumption $e_i \notin f(M')$ implies $e_i \notin f_Q(M'_Q)$: suppose that $e_i = \sum_{i' \leq n'} \kappa_{i'} f_{i'}$ where $(\kappa_{i'}) \in \cO_Q(n')$ and $f_{i'} := f(e_{i'})$ are the images of the basis vectors of $M'$. By \refle{particular}, we have $f_{i'} = \epsilon_{i'} e_i$ for all $i'$, with some $\epsilon_{i'} \in \cO_Q$, $|\epsilon_{i'}|=1$. But $f_{i'}$ also lies in $M$ (as opposed to $M_Q$). Thus, $\epsilon_{i'}$ must lie in $\cO$, that is, $e_i \in f(M')$. Therefore, to prove the claim we may assume $K$ is a field. 

Now, by \refle{particular}, $e_i$ is not a non-trivial $\cO$-linear combination of other elements of $M$. As $e_i \notin f(M')$, \refpr{cokernel} implies
\eqn \mylabel{eqn_pi-1}
\pi^{-1}(\pi(e_i)) = \{ e_i \}.
\xeqn

Fix a section $\sigma: M'' \r M$ of $\pi$, which exists by the assumption that $M''$ be projective. We obtain $\sigma(\pi(e_i)) = e_i$. Hence,
$$0 = \sigma (0_{M''}) = \sigma (\pi(t m'_1)) = \sum_{e_i \notin f(M')} t \lambda_i \sigma (\pi(e_i)) = \sum_{e_i \notin f(M')} t \lambda_i e_i,$$
so that $\lambda_i = 0$. The claim is shown.

By the claim, $f$ induces a bijection $\tilde f: M' = \cO(n') \r \cO(\tilde n')$, which gives rise to a bijection $K^{n'} \r K^{\tilde n'}$. This shows $\tilde n' = n'$. We conclude that the basis vectors $e_i \in M'$ get mapped under $f$ to $\epsilon_i e_{J(i)}$ where $\epsilon_i \in E$ and $J : \{1, \dots, n'\} \r \{1, \dots, n\}$ is an injective set map. In fact, suppose $\tilde f^{-1}(e_i) = \sum_{j \in J} \lambda_{ij} e_j$ with $(\lambda_{ij}) \in \cO(J)$ with all $\lambda_{ij} \neq 0$. Equivalently, $\sum \lambda_{ij} \tilde f(e_j) = e_i$. 
Therefore, by \refle{particular} (applied with $Q$ instead of $K$), $\tilde f_Q(e_j) \in E_Q \cdot e_i$ for all $j$, where $E_Q = \{ q \in Q, |q| = 1\}$. Since $\tilde f$ and therefore, by \refle{morphismscO}, $\tilde f_Q$ is injective, this implies that only one summand appears in this sum, i.e., 
$\tilde f(e_j) = \lambda_{ij}^{-1} e_i$ for some $j \in J$. A priori, $\lambda_{ij}^{-1}$ only lies in $Q$, but $\tilde f(e_j) \in \cO(n')$ shows that $\epsilon_i := \lambda_{ij}^{-1} \in \cO$, hence in $E$.

By \refas{basic}, $\epsilon_i \in E$ is a unit in $K$. We can therefore define
%
%
%
%
$\phi': \cO(n') \r M'$ by mapping the basis vectors $e_{i}$ of $\cO(n')$ (which correspond, in the above notation, to the basis vectors $e_{J(i)}$ of $M$) to $\epsilon_i^{-1} e_i$. Also, let $\phi'': \cO(n-n') \subset M \r M''$ be the map which sends the remaining basis vectors $e_{j'}$ for $j' \notin \im J$ to $\pi(e_{j'})$. Put
$$\phi := \phi' \sqcup \phi'': M = \cO(n) = \cO(n') \sqcup \cO(n-n') \r M' \sqcup M''.$$ 
Both $\phi'$ and $\phi''$ are onto, hence so is $\phi$. This follows from the construction of coproducts of modules over generalized rings \cite[4.6.15]{Durov}. (Also see \cite[10.4.7]{Durov} for an explicit description of the coproduct for modules over archimedean valuation rings.) Alternatively, the surjective maps $\phi'$ and $\phi''$ are epimorphisms of $\cO$-modules. Hence their coproduct $\phi$ is an epimorphism. As $M' \sqcup M''$ is projective, $\phi$ has a section, so it is also surjective. 
The map $\phi$ is injective, as can be seen by checking the definition or using \refle{morphismscO}\refit{l2c} $\Rightarrow$ \refit{l2b}. Hence $\phi$ is an isomorphism.

We finally show the unicity of $\phi$ or, in other words, that there are no non-trivial automorphism of cofiber sequences
$$0 \r M' \rightarrowtail M \twoheadrightarrow M'' \r 0.$$
Suppose $\tilde \phi$ is another isomorphism fitting into \refeq{short.exact}. We replace $\phi$ by $\tilde \phi \phi^{-1}$ and $\tilde \phi$ by $\id_M$ and assume $f$ is the standard inclusion $M' \r M = M' \sqcup M''$ and $\pi$ is the standard projection onto $M''$. Applying the base change functor $(-)_Q$ (see above), we may assume that $K$ is a field. Then $M''_K$ is a free $K$-module, so the endomorphism $\phi_K: M_K \r M_K$ is given by a matrix
$$B = \left ( \begin{array}{cc}
\Id_{M'} & A \\
0 & \Id_{M''}
\end{array}
\right ),$$
where $A$ is the matrix corresponding to the map $M''_K \r M'_K$ (of free $K$-modules). On the other hand, $\phi$ is a map of free $\cO$-modules, so every column in $B$ is in $\cO(n)$. This forces $A = 0$, so that $\phi = \id_M$.
\xpf

\theo \mylabel{theo_Waldhausen}
The category $(\Free(\cO), \cofi, \weq)$ defined in \ref{defi_Waldhausen} is a Waldhausen category. 
\xtheo

\pf

The only non-trivial thing to show is the stability of cofibrations under cobase-change. By \refpr{morphismscO2}, a cofibration sequence $M' \stackrel{\iota}\cofi M \stackrel \pi \fib M''$ in $\Free(\cO)$ is isomorphic to $M' \cofi M' \sqcup M'' \fib M''$. Hence, given any map $f: M' \r \tilde M'$, the pushout of $\iota$ along $f$, $\tilde M' \r \tilde M' \sqcup_{M'} M$ is isomorphic to $\tilde M' \r \tilde M' \sqcup M''$ which is a monomorphism with cokernel $M''$. 
\xpf

\rema \mylabel{rema_split}
Mahanta uses split monomorphisms as cofibrations in the category of finitely generated modules over a fixed $\Fone$-algebra (i.e., pointed monoid) to define $G$- (a.k.a. $K'$-)theory of such algebras \cite{Mahanta:G-Theory}. In $\Free(\cO)$, we have seen that all cofibrations are split, but not conversely: the cokernel of the split monomorphism $\varphi: \Zinf(1) \r \Zinf(2)$, $e_1 \mapsto \frac {e_1}2 + \frac {e_2}2$ is not free. This follows either from \refpr{morphismscO2} or by an explicit computation, using \refpr{cokernel}. Indeed, two elements $x_i e_1 + y_i e_2 \in \Zinf(2)$ ($i=1,2$) are identified in $\coker \varphi$ iff $|y_1 - x_1| = |y_2 - x_2| < 1$. On $\coker \varphi$, multiplication with $1/2$ is therefore not injective. Thus $\coker \varphi$ is not a submodule of a free $\Zinf$-module, in particular it is not projective.
\xrema

\subsection{$K$-theory}
In this subsection, we compute the $K$-theory of the generalized valuation ring $\cO$ (\refde{cO}) or, more precisely, of the category of free $\cO$-modules. By \refth{projfree}, every projective $\cO$-module is free, provided that the norm is archimedean.

We define the $K$-theory using Waldhausen's $S_\bullet$-construction, which has the advantage of being immediately applicable (\refth{Waldhausen}). Other constructions, such as Quillen's $Q$-construction can also be applied (slightly modified, since $\cO$-modules do not form an exact category). The resulting $K$-groups do not depend on the choice of the construction. 

Recall the definition of $K$-theory of a Waldhausen category $\cC$ (see e.g. \cite[Section IV.8]{Weibel:Kbook} for more details). We always assume that the weak equivalences of $\cC$ are its isomorphisms. The category $S_n \cC$ consists of diagrams
\eqn \mylabel{eqn_objS}
\xymatrix{
0 = A_{00} \ar@{>->}[r] & A_{01} \ar@{>->}[r] \ar@{->>}[d] & A_{02} \ar@{>->}[r] \ar@{->>}[d] & \dots \ar@{>->}[r] & A_{0n} \ar@{->>}[d] \\
 & 0 = A_{11} \ar@{>->}[r] & A_{12} \ar@{>->}[r] \ar@{->>}[d] & \dots \ar@{>->}[r] & A_{1n} \ar@{->>}[d] \\
 & & 0 = A_{22} \ar@{>->}[r] & \dots \ar@{>->}[r] & A_{2n} \ar@{->>}[d] \\
 & & & \ddots & \vdots \ar@{->>}[d] \\\
 & & & & A_{n-1,n} \\
}
\xeqn
such that $A_{i,j} \cofi A_{i,k} \fib A_{j,k}$ is a cofibration sequence. Varying $n$ yields a simplicial category $S_\bullet \cC$. The subcategory of 
isomorphisms is denoted $w S_\bullet \cC$. Applying the classifying space construction of a category yields a pointed bisimplicial set $S(\cC)_{n, m} := B_m w S_n \cC$. For example, $S(\cC)_{n,0} = \text{Obj} (S_n \cC)$. 
The $K$-theory of $\cC$ is defined as
$$K_i(\cC) := \pi_{i+1} d (B_* w S_\bullet \cC),$$
where $d(-)$ is the diagonal of a bisimplical set. 

By \refth{Waldhausen}, we are ready to define the \emph{algebraic $K$-theory} of $\cO$. More precisely, we consider the Waldhausen category of (finitely generated) free $\cO$-modules, which is the same as projective $\cO$-modules in all cases of interest by \refth{projfree}.

\defi \mylabel{defi_K}
\eqnarr
K_i (\cO) & := & K_i (\Free(\cO)) = \pi_{i+1} (d B w S_\bullet \Free(\cO)), \ \ i \geq 0.
\xeqnarr 
\xdefi

\lemm \mylabel{lemm_functoriality}
Given two normed domains and a ring homomorphism $f: K \r K'$ between them satisfying $|f(x)| = |x|$ (so that $f$ restricts to a map $f: \cO \r \cO'$), the functor $f^*: \Free(\cO) \r \Free(\cO')$, $M \mapsto M \t_{\cO} \cO'$ is (Waldhausen-)exact and therefore induces a functorial map
$$f^* : K_i (\cO) \r K_i (\cO').$$ 
\xlemm
\pf
As pointed out at p. \pageref{colimits}, $f^*: \Mod(\cO) \r \Mod(\cO')$ preserves cokernels. Secondly, tensoring with $\cO'$ preserves cofibrations since a map $M \r M'$ of free (or projective) $\cO$-modules is a monomorphism iff $M_Q \r M'_Q$ is one (where $Q$ is the quotient field of $K$, \refle{morphismscO}) and the statement is true for $Q$-modules: the map $Q \r Q'$ is injective since $|f(1)| = |1| = 1$ and therefore flat. 
\xpf

The group $K_0(\cO)$ is the free abelian group generated by the isomorphisms classes of free $\cO$-modules modulo the relations
$$[\cO(n') \sqcup \cO(n'')] = [\cO(n')] + [\cO(n'')].$$
Indeed, any cofiber sequence satisfies additivity of the ranks of the involved free modules, as one sees by tensoring the sequence with the quotient field $Q$ of $K$. Therefore, $K_0(\cO) = \Z$. 

We now turn to higher $K$-theory of $\cO$. Recall that $E : = \{x \in \cO, |x| = 1 \}$ is the subgroup of norm one elements. Let us write $\GL_n (\cO) := \Aut_{\cO}(\cO(n))$. According to \refpr{morphismscO2}, 
\eqn \mylabel{eqn_GLn}
\GL_n(\cO) = E \wreath S_n = E^n \rtimes S_n,
\xeqn 
where the symmetric group $S_n$ acts on $E^n$ by permutations. For $E = \mu_2 = \{ \pm 1 \}$, this group is known as the \emph{hyperoctahedral group}. As usual, we write
$$\GL(\cO) := \varinjlim_n \GL_n(\cO)$$ 
for the infinite linear group, where the transition maps are induced by $\GL_n(\cO(n) \ni f \mapsto f\sqcup \id_{\cO}$. For any group $G$, let $G_\ab = G / [G, G]$ be its abelianization. 
We write $\pist{i}(-)$ for the stable homotopy groups of a space and abbreviate $\pist i := \pist i(S^0)$.

\theo \mylabel{theo_Khigher}
Let $\cO$ be a generalized valuation ring as defined in \ref{defi_cO}. Then for $i\geq 0$, there is an isomorphism
$$K_i(\cO) \cong \pist i (BE_+, *),$$
where the right hand side denotes the $i$-th stable homotopy group of the classifying space of $E$ (viewed as a discrete group), with a disjoint base point $*$. For a map $f$ as in \refle{functoriality}, this isomorphism identifies $f^*$ in $K$-theory with the map on stable homotopy groups induced by $E(\mathcal O) \r E(\mathcal O')$. 

For $i=1$, $2$ we get
\eqnarra
K_1(\cO) & = & \GL(\cO)_\ab = E \x \Z/2 \nonumber \\
K_2 (\cO) & = & \varinjlim_n \H_2([\GL_n(\cO), \GL_n(\cO)], \Z) \mylabel{eqn_K2}
\xeqnarra 
where the right hand side in \refeq{K2} is group homology with $\Z$-coefficients.
\xtheo

Before proving the theorem, we first discuss our main example, when $\cO$ comes from an infinite place of a number field, as in \refex{1}. Then, we prove a preliminary lemma.

\exam \mylabel{exam_Atiyah}
Let us consider a number field $F$ with the norm induced by some complex embedding $\sigma \in \Sigma_F$ (see p.\ \pageref{notation} for notation). The torsion subgroup $E_\text{tor}$ of $E := \{x \in F^\x, |x| = 1\}$ agrees with the finite group $\mu_F$ of roots of unity. The exact localization sequence involving all finite primes of $\OF$, 
$$1 \r \cO_F^\x \r F^\x \r L := \ker (\oplus_{\pp < \infty} \Z \r \text{cl}(F)) \r 0,$$
shows $F^\x / \mu_F \cong \cO_F^\x / \mu_F \oplus L$. Hence it is free abelian by Dirichlet's unit theorem. Thus 
$$E \subset \mu_F \oplus \Z^{r_1 + r_2 - 1} \oplus L,$$ 
where $r_1$ and $r_2$ are the numbers of real and pairs of complex embeddings. Therefore, $E = \mu_F \oplus \Z^S$, where $S := \rk E$ is at most countably infinite. Of course, $E = \{ \pm 1 \}$ whenever $\sigma$ is a real embedding, but also, for example, for any complex embedding of $F = \Q[\sqrt[3]{2}]$. For $F = \Q[\sqrt{-1}]$, $E$ is the (countably) infinitely generated group of pythagorean triples \cite{Eckert} (see also \cite{ZanardoZannier} for a description of the group structure of pythagorean triples in more general number fields).  

The group $\mu_F$ is cyclic of order $w$, so the long exact sequence of group homology,
$$\H_i(\mu_F, \Z) \stackrel{\cdot n}\lr \H_i(\mu_F, \Z) \r \H_i(\mu_F, \Z/n) \r \H_{i-1}(\mu_F, \Z),$$
together with the Atiyah-Hirzebruch spectral sequence
$$\H_p (\mu_F, \pist q) = \H_p (B \mu_F, \pist q) \Rightarrow \pist{p+q} (B \mu_F) = \pist{p+q} ((B \mu_F)_+, *)$$
yield at least for small $p$ and $q$ explicit bounds on $\pist{p+q}((B \mu_F)_+, *)$: the $E^2$-page reads
$$\begin{array}{c|cccc}
q \uparrow \\
2 & \pist 2 = \Z/2 & \Z/w' & \Z/w'\\
1 & \pist 1 = \Z /2 & \mu_F / 2 = \Z/w' & \Z/w'\\
0 & \Z & \mu_F = \Z/w & 0 \\
\hline
& 0 & 1 & 2 & p \r
  \end{array}
$$
where $w'=(2,w)$. In general, $\pist{p+q}((B \mu_F)_+, *)$ is finite for $p+q>0$. 
For $i > 0$,
\eqnarr
K_i(\OF_\sigma) & = &  K_i(\OF_{(\sigma)}) \\
& = & \pist i (B (\mu_F \oplus \Z^{\oplus S})_+, *) \\
& =&  \pist i ((B\mu_F)_+ \vee \bigvee_S S^1, *) \\
&=& \pist i (B\mu_F) \oplus \bigoplus_S \pist{i-1}.
\xeqnarr
In particular 
\eqnarr
K_1(\OF_{(\sigma)}) & = & \Z/2 \oplus \mu_F \oplus \Z^{\oplus S}, \\ 
K_2(\OF_{(\sigma)}) & = &  G \oplus (\Z/2)^{\oplus S},
\xeqnarr 
where $G$ is a finite (abelian) group which is filtered by a filtration whose graded pieces are subquotients of $\Z/2$ and $\Z/w'$. 
(Determining $G$ would require studying the differentials of the spectral sequence).  
\xexam

\lemm \mylabel{lemm_GLOab}
The map
$$\GL(\cO)_\ab \r E \x \Z/2, (\epsilon, \sigma) \mapsto (\prod_{i=1}^\infty \epsilon_i, \text{parity}(\sigma))$$
is an isomorphism. Here the representation of elements of $\GL(\cO)$ is as in \refeq{GLn}. The group $[\GL(\cO), \GL(\cO)]$ is perfect.
\xlemm

\pf
For $i \geq 1$ and $\epsilon \in E$, let $\epsilon_i = (1, \dots, 1, \epsilon, 1, \dots) \in E \x E \x \dots$ be the vector with $\epsilon$ at the $i$-th spot. Let $\sigma_i = (i, i+1) \in S_n$ be the permutation swapping the $i$-th and $i+1$-st letter. The $\epsilon_i$ and $\sigma_i$, for $i \geq 1$ and $\epsilon \in E$, generate $G := \GL(\cO)$ as we have seen in the proof of \refpr{morphismscO2}. In $G$, we have relations $\sigma_i \sigma_{i+1} \sigma_i = \sigma_{i+1} \sigma_i \sigma_{i+1}$, which implies 
$\sigma_i = \sigma_{i+1}$ in $G_\ab$. Moreover, in $G$ we have the relation $\epsilon_i \sigma_i = \sigma_{i+1} \epsilon_{i+1}$, so that we get $\epsilon_i = \epsilon_{i+1}$ in $G_\ab$. This shows the first claim. 

The perfectness of $[\GL(\cO), \GL(\cO)]$ is a special case of \cite[Prop. 3]{Weibel:KtheoryAzumaya}, for example. 
Alternatively, the above implies that $H := [\Aut(\cO(n)), \Aut(\cO(n))]$ is given by $H = L \rtimes A_n$, where the alternating group $A_n$ acts on $L := \ker (\prod_{i=1}^n E \r E, (\epsilon^1, \dots, \epsilon^n) \mapsto \prod \epsilon^i) (\cong E^{n-1})$ by restricting the $S_n$-action on $E^n$. Now, the perfectness of $A_n$ for $n \geq 5$ and a simple explicit computation shows $H_\ab=1$ for $n \geq 5$. 
\xpf

We now prove \refth{Khigher}. This theorem is actually an immediate consequence of \refpr{morphismscO2}, together with well-known facts about $K$-theory of $G$-sets, where $G$ is some group \cite[Ex. IV.8.9]{Weibel:Kbook}. For example, the $K$-theory of the Waldhausen category of finite pointed sets (which would correspond to the impossible case $E=1$) is
$$K_i(\Fone) := K_i ((\text{finite pointed sets}, \text{ injections, bijections})) = \pist i,$$
the stable homotopy groups of spheres. More generally, for some (discrete) group $G$, the $K$-theory of the category $\Free(G)$ of finitely generated (i.e., only finitely many orbits) pointed $G$-sets on which the $G$-action is fixed-point free, together with bijections as weak equivalences and injections as cofibrations, is known to be the stable homotopy group of $(BG)_+$. By \refpr{morphismscO2}, the canonical functor
$$\Free (E) \r \Free (\cO), (E^X) \sqcup \{ * \} \mapsto \cO(X)$$
induces an equivalence of the categories of cofibrations and therefore an isomorphism of $K$-theory. For the convenience of the reader, we recall the necessary arguments, which also includes showing that other definitions of higher $K$-theory (of free $\cO$-modules) yield the same $K$-groups.  

\pf
Let $Q \Free(\cO)$ be Quillen's $Q$-construction, i.e., the category whose objects are the ones of $\Free(\cO)$ and 
$$\Hom_{Q \Free(\cO)} (A, B) := \{ A \fibl A' \cofi B \} / \sim,$$ 
where two such roofs are identified if there is an isomorphism between them which is the identity on $A$ and $B$. It forms a category whose composition is given by the composite roof defined by the cartesian diagram
$$\xymatrix{
& & A'' := A'\x_B B' \ar@{->>}[dl] \ar@{>->}[dr] \\
& A' \ar@{->>}[dl] \ar@{>->}[dr] & & B' \ar@{>->}[dr] \ar@{->>}[dl] \\
A & & B & & C.
}$$
Here, we use that $A''$ exists (in $\Free(\cO)$) since it is the kernel of the composite $B' \fib B \fib B / A'$,  
which is split by \refpr{morphismscO2}. The subcategory $S := \text{Iso} (\Free(\cO))$ of $\Free(\cO)$ consisting of isomorphisms only is a monoidal category under the coproduct. Hence $S^{-1} S$ is defined. We claim
$$\Omega B Q \Free(\cO) = B (S^{-1} S).$$
Indeed, the proof of \cite[Theorem IV.7.1]{Weibel:Kbook} carries over: the extension category $\mathcal E \Free(\cO)$ is defined as in \lcs and comes with a functor $t: \mathcal E \Free(\cO) \r Q \Free(\cO), (A \cofi B \fib C) \mapsto C$. The fiber $\mathcal E_C := t^{-1} C$ ($C \in \Free(\cO)$) consists of sequences $A \cofi B \fib C$. The functor
$$\phi: S \r \mathcal E_C, \ \  A \mapsto A \cofi A \sqcup C \fib C$$
induces a homotopy equivalence $B(S^{-1} S) \r B(S^{-1} \mathcal E_C)$ in the classical case of an exact category (instead of $\Free(\cO)$). In our situation, $\phi$ is an equivalence of categories since any extension in $\Free(\cO)$ splits \emph{uniquely} (\refpr{morphismscO2}). Thus  \cite[Theorem IV.4.10]{Weibel:Kbook} gives
$$B Q \Free(\cO) = K_0(S) \x B \GL(\cO)^+,$$ 
where the right hand side is the $+$-construction with respect to the perfect normal subgroup $[\GL(\cO), \GL(\cO)]$ (\refle{GLOab}). In the same vein, Waldhausen's comparison of the $Q$-construction and his $S_\bullet$-construction carries over: $d (B w S_\bullet \Free(\cO))$ is weakly equivalent to $B Q \Free(\cO)$.  

Finally, by the Barratt-Priddy theorem (see e.g. \cite[Th. 3.6]{Segal:Categories})  
$$\pi_i (B \GL(\cO)^+) \cong \pist i (B E_+, *).$$

The identification of the low-degree $K$-groups is the standard calculation of the $S^{-1}S$-construction \cite[IV.4.8.1, IV.4.10]{Weibel:Kbook}. 
\xpf

\rema
The calculation of $K_1(\cO)$ could also be done using the description of $K_1$ of a Waldhausen category due to Muro and Tonks \cite{MuroTonks}. 
\xrema

\rema \mylabel{rema_cofinality}
Recall that for an (ordinary) ring $R$ the following two properties of an $R$-module $M$ are equivalent: (i) it is projective, (ii) there is another projective module $M'$ such that $M \sqcup M'$ is free. I have not been able to show the corresponding statement for projective $\cO$-modules. For example, for a projector $p: \cO(n) \r \cO(n)$ with $M = \im p$, it is \emph{not} true that the canonical map
$$\phi: M \sqcup \ker p \r \cO(n)$$
is an isomorphism of $\cO$-modules: for $n=2$ and the projector $p$ given by the matrix 
$$\left ( \begin{array}{cc} 1/2 & 1/2 \\ 1/2 & 1/2 \end{array} \right ),$$ 
$\ker p$ is the free $\cO$-module of rank $1$, generated by $(e_1 - e_2) /2 \in \cO(2)$. In this case, $\phi$ induces an isomorphism of $M \sqcup \ker p$ with the free $\cO$-module of rank $2$ generated by $(e_1 \pm e_2)/2$, but not with $\cO(2) = (e_1, e_2)$. 
The analogous statement of \refpr{morphismscO2} for cofibrations of projective $\cO$-modules, as well as the computation of $K_i(\Proj(\cO))$ for $i > 0$ (using Waldhausen's cofinality theorem) would carry over verbatim if the above statement about projective $\cO$-modules holds. However, the distinction between projective and free modules is only relevant for non-archimedean valuations, by \refth{projfree}. 
\xrema

\section{The residue field at infinity}
We finish this work by noting two differences (as far as $K$-theory is concerned) to the case of classical rings, namely the $K$-theory of the residue ``field'' at infinity, and the behavior with respect to completion. For simplicity, we restrict our attention to the case $F = \Q$. 

Let $p < \infty$ be a (rational) prime with residue field $\Fp$. There is a long exact sequence
$$K_n (\Fp) \r  K_n (\Z_{(p)}) \r K_n(\Q) \stackrel{\delta} \r K_{n-1} (\Fp)$$
which stems from the fact that $\Z_{(p)}$ (the localization of $\Z$ at the prime ideal $(p)$) is a Noetherian regular local ring of dimension one. Moreover, for $n=1$ the map $\delta$ is the $p$-adic valuation $v_p : \Q^\x \r \Z$. The situation is less formidable at the infinite places, as we will now see. The (generalized) valuation ring $\Zinff$ (\refde{cO}) is \emph{not} Noetherian: ascending chains of ideals need not terminate. Indeed, consider a finitely generated ideal $I = (m_1, \cdots, m_n) \subset \Zinff$. Then $|I| = \{ |m|, m \in I \} 
= [0, \max_i |m_i|] \cap |\Zinff|$. In particular, an ideal of the form $\{x \in \Zinff, |x| < \lambda \}$, 
$\lambda \leq 1$ is not finitely generated, since $|\Zinff|$ is dense in $[0,1]$. 
This should be compared with the well-known fact 
that the valuation ring of a non-archimedian field is noetherian iff the field is trivially or discretely valued. 

\defi \cite[4.8.13]{Durov}.
Put $\Finf := \Zinff / \widetilde {\Zinff}$, where $\widetilde {\Zinff}$ is the submonad given by
$$\widetilde {\Zinff}(n) = \{ x \in \Q^n, |x| < 1\}.$$
\xdefi
We refer to \lcs for the general definition of strict quotients of generalized rings by appropriate relations. For us, it is enough to note that every element of $\Zinff(n)$ is uniquely represented by $z = \sum_{i \in I} \lambda_i \epsilon_i e_i$, where $I \subset \{1, \dots, n\}$, $0 < \lambda_i \leq 1$, $\sum \lambda_i \leq 1$, $\epsilon_i \in E_{\Zinff} = \{ \pm 1 \}$, and $e_i$ is the standard basis vector. Two elements $z, z' \in \Zinff(n)$ get identified in $\Finf(n)$ (Notation: $z \equiv z'$) iff 
\eqn
\mylabel{eqn_ident1} |z| < 1 \text{ and } |z'|<1
\xeqn
or
\eqn
\mylabel{eqn_ident2} |z| = |z'|=1,\ I_z=I_{z'}, \text { and } \epsilon_{i, z}=\epsilon_{i, z'} \text{ for all }i \in I_z.
\xeqn
That is, as a set $\Finf(n)$ consists of the faces of the $n$-dimensional octahedron. Again, $0$ is the initial and terminal $\Finf$-module, so we can speak about (co)kernels.

As usual, we put 
$$K_0(\Finf) := \left (\bigoplus_{M \in \Free(\Finf) / Iso} \Z  \right) \left / [M] = [M'] + [M''] \right .,$$
with a relation for each monomorphism $M' \r M$ in $\Free(\Finf)$ such that its cokernel $M''$ (computed in $\Mod(\Finf)$) lies in $\Free(\Finf)$. Similarly, we define $K_0^{\Proj}(\Finf)$ using projective $\Finf$-modules. Using the above, one sees that $\Finf$ is not finitely presented as a $\Zinff$-module.  Thus, one should not expect a natural map $i_*: K_0(\Finf) \r K_0(\Zinff)$. Actually, $K$-theory of $\Finf$-modules behaves badly in the sense of the following proposition:

\prop \mylabel{prop_K0Finf}
$K_0^\Proj(\Finf) = 0$, $K_0(\Finf) = \Z$.
In particular, there is no exact \emph{localization sequence} (regardless of the maps involved)
$$K_1(\Zinff) = \Z/2 \x \{ \pm 1 \} \r K_1(\Q) = \Q^\x \r K_0(\Finf) \r K_0(\Zinff) = \Z \r K_0(\Q) = \Z,$$
or similarly with $K_0^\Proj(\Finf)$ instead.  
\xprop
\pf
We first show that any projective $\Finf$-module $M$ which is generated by $n$ elements contains $\Finf$ as a submodule, such that the cokernel is a projective $\Finf$-module generated by $n-1$ elements. This implies that $K^\Proj_0(\Finf)$ is generated by $[\Finf]$ (which is obvious for $K_0(\Finf)$).

The projective module $M$ is specified by a projector $\pi : \Finf(n) \r \Finf(n)$ with $M = \pi(\Finf(n))$. Let $a_i := \pi(e_i) \in \Finf(n)$. We pick $a_{ij} \in [-1, 1] \subset \RR$ such that $a_i \equiv \sum_{j \in J_i} a_{ij} e_j$ with $a_{ij} \neq 0$ for all $j \in J_i$. Set $A := (a_{ij}) \in \RR^{n\x n}$. We may assume that the number $n$ of generators of $M$ is minimal, i.e., there is no surjection $p': \Finf(n') \r M$ with $n' < n$. 
Indeed, if there is such a surjection, it has a section $\sigma'$ since $M$ is projective, and $\pi' := \sigma' p'$ would again be a projector.

The minimality of $n$ implies that $a_i \nequiv a_j$ for all $i \neq j$. Otherwise, the restriction of $\pi$ to $\Finf (n \backslash \{ i \}) \subset \Finf(n)$ would be surjective. Similarly, the minimality implies $a_i \nequiv 0 \in \Finf(n)$ for all $i$. 
Also, put $B = (b_{ij}) := A^2 \in \RR^{n \x n}$. Using $(b_{ij})_j \equiv \pi(a_i) \equiv a_i \nequiv 0 \in \Finf(n)$, we obtain $\sum_{j} |b_{ij}| = 1$ and $\sum_j |a_{ij}| = 1$ by \refeq{ident1}. 

The minimality of $n$ implies $i \in J_i$ or equivalently, $a_{ii} \neq 0$: otherwise $a_i \equiv \pi (a_i) \equiv \sum_{j \in J_i \backslash \{ i\}} a_{ij} a_j$ would be an $\Finf$-linear combination of the remaining columns of $A$.

For every $i \leq n$,
\eqnarr
1 & = & \sum_j |b_{ij}| = \sum_j | \sum_k a_{ik} a_{kj}| \\
& \leq & \sum_j \sum_k |a_{ik}| |a_{kj}|  = \sum_k |a_{ik}| \underbrace{\left (\sum_j |a_{kj}|\right)}_{=1} \\
& = & 1,
\xeqnarr
so equality holds. In particular, the terms $\sgn (a_{ik} a_{kj})$ are either all (for arbitrary $i, j, k \leq n$) non-negative or non-positive. Picking $k = j := i$, we see that they are non-negative, since $\sgn (a_{ii}^2) > 0$, for $a_{ii} \neq 0$. 

Let $I^> := \{i, a_{ii} > 0\}$ and likewise with $I^<$. Then $I^> \sqcup I^- = \{ 1, \dots, n \}$. Moreover, for $i \in I^>$ and $j \in I^<$, $a_{ii} a_{ij} \geq 0$ and $a_{ij} a_{jj} \geq 0$ imply $a_{ij} = 0$. In other words, the matrix $A$ decomposes as a direct sum matrix $A^> \sqcup A^<$, where $A^>$ and $A^<$ are the submatrices of $A$ consisting of the rows and columns with indices in $I^>$ and $I^<$, respectively. We may therefore assume $A = A^>$, say. For $i (\in I^>)$, and any $j$, $a_{ii} a_{ij} \geq 0$ implies $a_{ij} \geq 0$, i.e., the entries of $A$ are all non-negative. 
 
Fix some $i \leq n$. As $\pi$ is a projector, $a_i \equiv \pi(a_i)$, i.e.,
$$a_i \equiv \sum_{j \in J_i} a_{ij} e_j \equiv  \sum a_{ij} \pi(e_j) \equiv  \sum_{j \in J_i, k \in J_j} a_{ij} a_{jk} e_k \in \Finf(n).$$
By \refeq{ident1}, \refeq{ident2}, this implies $\sgn(a_{ik}) = \sgn (\sum_j a_{ij} a_{jk})$, which gives
\eqn \mylabel{eqn_J} 
J_i = \cup_{j \in J_i} J_j.
\xeqn 
Indeed, ``$\subset$'' is easy to see without using the non-negativity of the entries. Conversely, for $k \notin J_i$, $\sum_j a_{ij} a_{jk} = 0$. Since all $a_{**} \geq 0$, this implies $a_{jk} = 0$ for all $j \in J_i$, i.e., $k \notin \cup_{j \in J_i} J_j$.

Now, pick some $i \leq n$ such that $J_i$ is maximal, i.e., not contained in any other $J_j$, $i \neq j$. Then $i \notin J_j$ for any $i \neq j$ by \refeq{J}. In other words, the $i$-th row only contains a single non-zero entry. For simplicity of notation, we may suppose $i = 1$.

Consider the diagram
$$\xymatrix{
\Finf \ar@{>->}[r]^\iota \ar@{=}[d] & \Finf(n) \ar@{->>}[r]^\rho \ar@{->>}[d] & \Finf(n-1) \ar@{->>}[d] \\
\Finf \ar@{>->}[r] & M \ar@{->>}[r] & M'
}$$
where $\rho$ is the projection onto the last $n-1$ coordinates, $\iota$ is the injection in the first coordinate. The lower left-hand map is a monomorphism since the first row of $A$ is nonzero. Its cokernel $M'$ is the projective module determined by the matrix $(a_{ij})_{2 \leq i,j \leq n}$. This exact sequence shows that $K^\Proj_0(\Finf)$ is generated by $[\Finf]$. 

On the other hand, consider the projective $\Finf$-module $P$ defined by the projector 
$\left (\begin{array}{cc} 1/2 & 0 \\ 1/2 & 1 \end{array} \right )$ \cite[10.4.20]{Durov}. It consists of $5$ elements and can be visualized as 
$$P = 
\vcenter
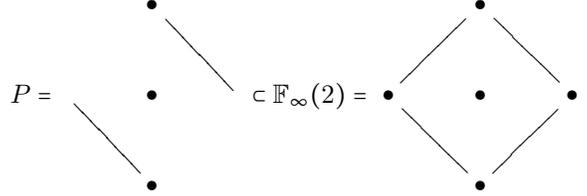
{\xymatrix{
& \bullet \ar@{-}[dr] & \\
\ar@{-}[dr] & \bullet &  \\
& \bullet &}} 
\subset
\Finf(2) = 
\vcenter{\xymatrix{
& \bullet \ar@{-}[dr]\ar@{-}[dl] & \\
\bullet \ar@{-}[dr] & \bullet & \bullet \ar@{-}[dl]  \\
& \bullet &}}.
$$ 
The composition $\Finf \stackrel{(1/2, 1/2)}\lr 
\Finf(2) \twoheadrightarrow P$ is a monomorphism with cokernel $\Finf$. The pictured inclusion $P \r \Finf(2)$ has cokernel $\Finf$, spanned by $e_1$. This shows that $[\Finf(2)] = 2 [\Finf] = [P] + [\Finf] = 3 [\Finf]$. Hence $K^\Proj_0(\Finf) = 0$.

Finally, we have to show $K_0(\Finf) = \Z$. For this, consider a cofiber sequence
$$\Finf(n') \stackrel i \cofi \Finf(n) \stackrel p \fib \Finf(n'').$$
We have to show $n = n' + n''$. Pick a section $\sigma$ of $p$. The natural map $i \sqcup \sigma: \Finf(n') \sqcup \Finf(n'') \r \Finf(n)$ is injective, as one easily shows. Thus $n' + n'' \leq n$ for cardinality reasons. Conversely, for any basis vector $e_i \in \Finf(n) \backslash \im i$, $p^{-1} (p(e_i)) = \{ e_i \}$, as one shows in the same way as for $\Zinf$-modules, cf.\ \refeq{pi-1}. Thus $\sigma(p(e_i)) = e_i$, so there are at most $n''$ such basis vectors by the injectivity of $\sigma$. Moreover, at most $n'$ of the basis vectors $e_i$ of $\Finf(n)$ are in $\im i$ by the injectivity of $i$. This shows $n' + n'' \geq n$. 
\xpf

\rema 
For $p \leq \infty$, let $Fib$ be the homotopy fiber of $\Omega K(\Z_{(p)}) \r \Omega K(\Q)$ and $\widehat {Fib}$ the one of  $\Omega K(\Zp) \r \Omega K(\Qp)$. The localization sequence for $K$-theory shows in case $p<\infty$ that $Fib$ and $\widehat {Fib}$ are homotopy equivalent (and given by $K(\Fp)$). Here $\Omega$ is the loop space and $K(-)$ is a space (or spectrum) computing $K$-theory, for example the $S_\bullet$-construction. 
However, for $p=\infty$, we have
$$\xymatrix{
\pi_1(Fib) \ar[r] & K_1(\Zinff) \ar@{=}[d] \ar[r] & K_1(\Q)=\Q^\x \ar[d]^{\subsetneq} \ar[r] & \pi_0(Fib) \ar[r] & 0 \\ 
\pi_1(\widehat {Fib}) \ar[r] & \underbrace{K_1(\Zinf)}_{(\Z/2)^{\oplus 2}} \ar[r] & K_1(\RR)=\RR^\x \ar[r] & \pi_0(\widehat {Fib}) \ar[r] & 0,  
}$$
so that $\pi_0(Fib) \subsetneq \pi_0(\widehat {Fib})$. 
\xrema  

\bibliography{bib}
\end{document}